\documentclass[12pt]{amsart} 
\usepackage[T1]{fontenc}     
\usepackage[utf8]{inputenc}  
\usepackage[margin=2.8cm]{geometry}   
\usepackage{amssymb}         
\usepackage{xcolor}          

\newtheorem{thm}{Theorem}
\newtheorem{cor}[thm]{Corollary}
\newtheorem{lem}[thm]{Lemma}
\theoremstyle{remark}
\newtheorem{rem}[thm]{Remark}
\newtheorem{exams}[thm]{Examples}

\newcommand{\N}{\mathbb{N}}
\newcommand{\R}{\mathbb{R}}
\newcommand{\Z}{\mathbb{Z}}
\newcommand{\T}{\mathbb{T}}
\newcommand{\calO}{\mathcal{O}}

\newcommand{\wrt}{\,\mathrm{d}}
\newcommand{\afterbar}{\bar{\phantom{n}}}

\begin{document}

\title{Open mappings of locally compact groups}
\author[M. G. Cowling]{Michael G. Cowling}
\address{School of Mathematics and Statistics, University of New South Wales, UNSW Sydney NSW 2025, Australia}
\email{m.cowling@unsw.edu.au}

\author[K. H. Hofmann]{Karl H. Hofmann}
\address{Fachbereich Mathematik, Technische Universität Darmstadt, Schlossgartenstraße 7, D-64289 Darmstadt, Germany}
\email{hofmann@mathematik.tu-darmstadt.de}

\author[S. A. Morris]{Sidney A. Morris}
\address{School of Mathematical and Physical Sciences, La Trobe University, Bundoora VIC 3086, Australia}
\email{morris.sidney@gmail.com}

\thanks{We are grateful to Gábor Lukács for his valuable comments on earlier versions of the text.  
The first author was partially supported by the Australian Reasearch Council (DP220100285).}

\begin{abstract}
The aim of this note is to insert in the literature some easy but apparently not widely known facts about morphisms of locally compact groups, all of which are concerned with the openness of the morphism.
\end{abstract}

\maketitle

\section{Main result}
We are interested in \emph{morphisms of locally compact groups}, by which we mean continuous homomorphisms between locally compact Hausdorff topological groups, and in particular when these are \emph{open}.

\begin{exams}\label{exam:three}
To set the scene, we recall a few standard examples and counterexamples of such morphisms $\phi\colon G\to H$:
\begin{enumerate}
  \item[(a)]
Let $\Z_p$ be the compact additive group of $p$-adic integers and $\R$ be the additive group of real numbers. 
Let $D$ be the discrete subgroup $\{ (m,m) : m \in\Z\}$ of $G :=\Z_p\times \R$.
Finally, define $H := G/D$, the compact \emph{$p$-adic solenoid} $\T_p$ (see \cite[Example 1.28 and Exercise E1.11]{HM1}), and let $\phi\colon G\to H$ be the quotient morphism. 
Then $\phi$ is a surjective open morphism of locally compact groups mapping the identity component $G_0=\{0\}\times \R\cong \R$ 
onto the dense arc component properly contained in the compact connected
group $H=H_0$ (see \cite[E1.11 (iii) and (iv)]{HM1}).
  \item[(b)]
Let $H$ be any nondiscrete locally compact group and $G$ be $H_d$, the underlying group of $H$ with the discrete topology. Then the identity map $\phi\colon G\to H$ is a nonopen bijective morphism.
  \item[(c)]
Let $G$ be any locally compact noncompact abelian group, $H$ be its Bohr compactification, and $\phi\colon G\to H$ be the natural standard morphism (see \cite[Chapter 8, Part 7]{HM1}). 
Then $\phi$ is a nonopen injective morphism with dense image properly contained in $H$.
This is related to the previous example by Fourier duality.
\end{enumerate}
\end{exams}

We write $\calO(G)$ for the space of open subgroups (which are automatically closed) of $G$, and $G \leqslant H$ (or $H \geqslant G$) to indicate that $G$ is a subgroup of $H$.
We recall that if $U$ is a subset of a locally compact group $G$ that is symmetric (meaning that $U^{-1} = U$), then $U$ generates a subgroup $\langle U \rangle$ of $G$ equal setwise to $\bigcup_{n \in \N} U^n$.
If $U$ is open, then $\langle U \rangle$ is open and if $U$ is relatively compact, then $\langle U \rangle$ is $\sigma$-compact.

\begin{thm}\label{thm:good-or-bad}   
Let $\phi\colon G\to H$ be a morphism of locally compact  groups.
If $\phi(G_1)$ is open in $H$ for all $G_1$ in $\calO(G)$, then $\phi$ is an open mapping.
\end{thm}

\begin{rem}
If $\phi$ is open, then trivially $\phi(G_1)$ is open in $H$ for all open subgroups $G_1$ of $G$, and so the theorem may also be stated as an equivalence.
\end{rem}

\begin{proof}
Take a symmetric relatively compact open subset $U$ of $G$; then $G_1 := \langle U \rangle$ is open, compactly generated, and $\sigma$-compact.

To show that $\phi\colon G \to H$ is open, it is necessary and sufficient to show that the restricted surjective morphism $\phi|_{G_1}\colon G_1 \to \phi(G_1)$ is open.
Indeed, if $\phi$ is open, then $\phi|_{G_1}$ is trivially open; conversely, every open set $U$ in $G$ is a union $\bigcup_{x \in U} x U_x$, where $U_x$ is an open subset of $G_1$; if $\phi|_{G_1}$ is open, then $\phi(x U_x)$ is open, and so $\phi(U)$ is open.

The Open Mapping Theorem for locally compact groups of \cite[Exercise EA1.21]{HM1} uses Baire category to show that, if $\phi\colon G \to H$ is a surjective morphism of locally compact groups and $G$ is $\sigma$-compact, then $\phi$ is open.
We apply this result to $\phi|_{G_1}$ and the theorem is proved.
\end{proof}

We recall that a subgroup $G_1$ of $G$ is open if and only if it contains the identity component $G_0$ of $G$ and $G_1/G_0$ is open in the totally disconnected group $G/G_0$.
Further, these open subgroups $G_1/G_0$ form a basis of the identity neighbourhoods of $G/G_0$, and 
\[
G_0 = \bigcap_{G_1 \in \calO(G)}  G_1.
\]
By a theorem of van Dantzig \cite{vD36}, whose proof is not very difficult, the $G_1/G_0$ that are both compact and open also form a basis of the identity neighbourhoods of $G/G_0$.

If $\phi\colon G\to H$ is a morphism, then $\phi(G_0)$ is connected so $\phi(G_0) \leqslant H_0$.
Hence $\phi$ induces a canonical morphism $\gamma\colon G/G_0 \to H/H_0$, given by $\gamma(g) = \phi(g)H_0/H_0$ for all $g \in G$.

\begin{cor} 
Assume that $\phi\colon G \to H$ is a morphism of locally compact groups and define $\gamma\colon G/G_0 \to H/H_0$ as above. 
Then the following are equivalent:
\begin{enumerate}
  \item[(1)] $\phi$ is open;
  \item[(2)] $\gamma\colon G/G_0 \to H/H_0$ is open and $H_0 = \bigcap_{G_1 \in \calO(G)} \phi(G_1)$.
\end{enumerate}
\end{cor}

\begin{proof}
Assume that (1) holds. 
For all $G_1 \in \calO(G)$, we have $\phi(G_1) \in \calO(H)$ and hence 
$H_0 \leqslant \phi(G_1)$.
Further, if $G_0\leqslant G_1\leqslant G$ and $G_1/G_0\in \calO(G/G_0)$, then $G_1\in \calO(G)$, and so
\[
\gamma(G_1/G_0)=\phi(G_1)H_0/H_0 = \phi(G_1)/H_0,
\]
and $\gamma(G_1/G_0)$ is open in $H/H_0$, that is, $\gamma$ is open.
Further, $\phi^{-1}(H_1) \in \calO(G)$ for all $H_1 \in \calO(H)$, and so
\[
H_0 \leqslant \bigcap_{G_1 \in \calO(G)} \phi(G_1) \leqslant \bigcap_{H_1 \in \calO(H)} H_1 = H_0.
\]
It follows that (2) holds.

Next, assume that $\Gamma$ is open and $H_0 \leqslant \phi(G_1)$ for all $G_1 \in \calO(G)$.
If $G_1 \in \calO(G)$, then $G_1/G_0 \in \calO(G/G_0)$, so $\gamma(G_1/G_0) \in \calO(H/H_0)$.
Since $H_0 \leqslant \phi(G_1)$, 
\[
\gamma(G_1/G_0) = \phi(G_1)H_0/H_0 = \phi(G_1)/H_0,
\]  
whence $\phi(G_1) \in \calO(G)$.
By Theorem \ref{thm:good-or-bad}, $\phi$ is open, as claimed.
\end{proof}

Examples show that ``simpler'' versions of this corollary fail.
On the one hand, if $G_0$ and $H_0$ are both open subgroups, then $\gamma$ is trivially open.
In this case, $\phi$ is open if and only if the restricted map $\phi|_{G_0}$ is open, and this holds if and only if $H_0 = \phi(G_0)$.
On the other hand, Example 1 (a) shows that there are open maps $\phi\colon G \to H$ for which $\phi(G_0) \neq H_0$.

Next, we consider where the chief difficulties of the theorem lie.
Assume that $\phi\colon G \to H$ is a morphism, and define $K := \ker(\phi)$ and $H_1 := \phi(G)\afterbar$.
Then $\phi$ factorises as the composition $\iota \circ \beta \circ \pi$, where $\pi$ is the canonical projection of $G$ onto $G/K$, $\beta\colon \phi(G) \to H_1$ is the injection of a dense subgroup of $H_1$ into $H_1$, and $\iota\colon H_1 \to H$ is the injection of the closed subgroup $H_1$ into $H$.

Since $\pi$ is open by definition, $\phi$ is open if and only if the canonical induced morphism $\psi\colon G/K \to H$ is open. 
Thus, in considering questions of openness of morphisms, there is no significant loss of generality in restricting our attention to the case where $\phi$ is injective and $\pi$ is trivial.

Likewise, if $H_1$ is not open, it is not possible that $\phi(U)$ can be open for any nonempty open subset $U$ of $G$, and there is no real loss of generality in restricting our attention to the case in which $H_1$ is open.
We may then assume that $H = H_1$, that is, that $\iota$ is trivial.

Finally, if $\iota$ and $\pi$ are trivial and $\phi$ sends open subgroups to open subgroups, then $\phi(G)$ is open and dense, and hence $\phi$ is surjective.
Thus the heart of the matter is what happens with groups $G$ and $H$ that are isomorphic as abstract groups, and the next corollary is actually an equivalent version of Theorem \ref{thm:good-or-bad}.

\begin{cor}    
Let $G$ and $H$ be locally compact groups, and let $\phi\colon G \to H$ be a continuous bijection of $G$ onto $H$. 
Assume that $\phi(G_1) \in \calO(H)$ for all $G_1 \in \calO(G)$.
Then $\phi$ is an isomorphism of locally compact groups, that is, $G\cong H$.
\end{cor}

\begin{proof}
By Theorem \ref{thm:good-or-bad}, the bijective morphism $\phi$ is open and thus is an isomorphism.
\end{proof}

\section{Another version}

In this section, we present a measure theoretic proof of Theorem \ref{thm:good-or-bad} and related results.
We begin with a simple lemma.

\begin{lem}\label{lem:good-or-bad}
Let $\phi\colon G \to H$ be a morphism of locally compact groups.
If $\phi$ is open, then $\phi(U)$ has nonzero measure in $H$ for all nonempty open subsets $U$ of $G$.
If $\phi$ is not open, then $\phi(V)$ has null measure in $H$ for all $\sigma$-compact subsets $V$ of $G$.
\end{lem}

\begin{proof}
If $U$ and $V$ are nonempty relatively compact (or $\sigma$-compact) open sets in $G$, then $V$ may be covered by finitely (or countably) many translates of $U$ and vice versa, and $\phi(V)$ may be covered by finitely (or countably) many translates of $\phi(U)$  and vice versa.
It follows that $\phi(U)$ has null measure if and only if $\phi(V)$ has null measure.

In light of this observation and the fact that the Haar measure of a nonempty open set is never $0$, it suffices to prove that if $U$ is a nonempty relatively compact open set in $G$ and $\phi(U)$ has positive Haar measure in $H$, then $\phi(U)$ is open in $H$.

Given any $z \in U$, $z^{-1}U$ is a relatively compact open neighbourhood of the identity in $G$, and by \cite[p.~18]{HR63}, it is possible to find a relatively compact open neighbourhood $V$ of the identity in $G$ that is symmetric ($V^{-1} = V$) such that $V^2 \subseteq z^{-1}U$.
It will suffice to show that $\phi(V^2)$ contains an open neighbourhood $W$ of the identity in $H$, for then $\phi(z)W$ is an open neighbourhood of $\phi(z)$ in $H$ that is contained in $\phi(zV^2)$, and hence in $\phi(U)$.

Since $V$ is relatively compact and $\phi$ is continuous, $\phi(V)$ is relatively compact in $H$ and has finite Haar measure, which is not $0$ by our initial observation, since the measure of $\phi(U)$ is not $0$.
Then the convolution on $H$ of the characteristic function of $\phi(V)$ with itself, given by the formula
\[
\chi_{\phi(V)} \ast \chi_{\phi(V)}(x)
= \int_H \chi_{\phi(V)}(y) \, \chi_{\phi(V)}(y^{-1}x) \wrt{y},
\]
is a continuous function on $H$ (see \cite[(20.16)]{HR63}) that does not vanish at the identity but vanishes outside $\phi(V^2)$, and $\{y \in H: \chi_{\phi(V)} \ast \chi_{\phi(V)}(y) \neq 0\}$ is the desired set $W$.
\end{proof}

Lemma \ref{lem:good-or-bad} has the following corollaries, the first of which we already know.

\begin{cor}\label{cor:good-or-bad-2}   
Let $\phi\colon G\to H$ be a morphism of locally compact  groups.
If $\phi(G_1)$ is open in $H$ for all compactly generated $G_1$ in $\calO(G)$, then $\phi$ is an open mapping.
\end{cor}

\begin{proof}
Assume that $\phi$ is not open, and take a nonempty relatively compact symmetric open subset $U$ of $G$. 
Then $U$ generates an open subgroup $\langle U \rangle$ of $G$, equal to $\bigcup_{n \in \N} U^n$. 
Since $U^n$ is a nonempty relatively compact open subset of $G$ for all positive integers $n$ and $\phi$ is not open, each $\phi(U^n)$ has measure $0$ by Lemma \ref{lem:good-or-bad}, whence $\phi(\langle U\rangle)$ has measure $0$, and so $\phi(\langle U\rangle)$ cannot be open.
\end{proof}

\begin{cor}    
Let $\phi\colon G \to H$ be a morphism of locally compact groups, and suppose that there is a compact subset $K$ of $G$ such that $\phi(K)$ has positive Haar measure in $H$.
Then $\phi$ is open.
\end{cor}

\section{Comments}
Our arguments and the proof of the Open Mapping Theorem require that we deal with a property of subsets of $G$ that is closed under taking countable unions, which is why Baire category or Haar measure play a role.

The reference \cite{HM2} provides an Open Mapping Theorem for almost connected pro-Lie groups
(Theorem 8.60, p.365) and an Alternative Open Mapping Theorem (Theorem 11.79, p. 517); earlier discussion of open mapping Theorems may be found in \cite{HM0a} and \cite{HM0b}.
This result continues to attract interest: see, for example, \cite{Ma23}.


\end{document}